
\magnification=1200
\catcode`\@=11
 
\def\hexnumber@#1{\ifnum#1<10 \number#1\else
 \ifnum#1=10 A\else\ifnum#1=11 B\else\ifnum#1=12 C\else
 \ifnum#1=13 D\else\ifnum#1=14 E\else\ifnum#1=15 F\fi\fi\fi\fi\fi\fi\fi}
 
\font\tenmsa=msam10
\font\sevenmsa=msam7
\font\fivemsa=msam5
\newfam\msafam
\textfont\msafam=\tenmsa  \scriptfont\msafam=\sevenmsa
  \scriptscriptfont\msafam=\fivemsa
\def\msa@{\hexnumber@\msafam}

\font\tenmsb=msbm10
\font\sevenmsb=msbm7
\font\fivemsb=msbm5
\newfam\msbfam
\textfont\msbfam=\tenmsb  \scriptfont\msbfam=\sevenmsb
  \scriptscriptfont\msbfam=\fivemsb
\def\msb@{\hexnumber@\msbfam}

\font\teneuf=eufm10
\font\seveneuf=eufm7
\font\fiveeuf=eufm5
\newfam\euffam
\textfont\euffam=\teneuf \scriptfont\euffam=\seveneuf
  \scriptscriptfont\euffam=\fiveeuf
\def\euf@{\hexnumber@\euffam}
\def\Frak{\relax\ifmmode\let\next\Frak@\else
 \def\next{\errmessage{Use \string\Frak\space only in math mode}}\fi\next}
\def\Frak@#1{{\Frak@@{#1}}}
\def\Frak@@#1{\fam\euffam#1}
 
\def\Bbb{\relax\ifmmode\let\next\Bbb@\else
 \def\next{\errmessage{Use \string\Bbb\space only in math mode}}\fi\next}
\def\Bbb@#1{{\Bbb@@{#1}}}
\def\Bbb@@#1{\fam\msbfam#1}
 
\font\tenscr=rsfs10 
\font\sevenscr=rsfs7 
\font\fivescr=rsfs5 
\skewchar\tenscr='177 \skewchar\sevenscr='177 \skewchar\fivescr='177
\newfam\scrfam
\textfont\scrfam=\tenscr \scriptfont\scrfam=\sevenscr
  \scriptscriptfont\scrfam=\fivescr
\def\Scr{\relax\ifmmode\let\next\Scr@\else
 \def\next{\errmessage{Use \string\Scr\space only in math mode}}\fi\next}
\def\Scr@#1{{\Scr@@{#1}}}
\def\Scr@@#1{\fam\scrfam#1}

\catcode`\@=12
 
\def\E{{\Bbb E}}
\def\L{{\Bbb L}}

\def\Q{{\Bbb Q}}

\def\C{{\Bbb C}}


\newcount\thmno \thmno=1
\long\def\theorem#1#2{\edef#1{Theorem~\the\thmno}\noindent
    {\bf Theorem \the\thmno.\enspace}{\it #2}\endgraf
    \penalty55\global\advance\thmno by 1}
\newcount\lemmano \lemmano=1
\long\def\lemma#1#2{\edef#1{Lemma~\the\lemmano}\noindent
    {\bf Lemma \the\lemmano.\enspace}{\it #2}\endgraf
    \penalty55\global\advance\lemmano by 1}
\newcount\propno \propno=1
\long\def\proposition#1#2{\edef#1{Proposition~\the\propno}\noindent
    {\bf Proposition \the\propno.\enspace}{\it #2}\endgraf
    \penalty55\global\advance\propno by 1}
\newcount\corolno \corolno=1
\long\def\corollary#1#2{\edef#1{Corollary~\the\corolno}\noindent
    {\bf Corollary \the\corolno.\enspace}{\it #2}\endgraf
    \penalty55\global\advance\corolno by 1}

\newcount\conjno \conjno=1
\long\def\conjecture#1#2{\edef#1{Conjecture~\the\conjno}\noindent
    {\bf Conjecture \the\conjno.\enspace}{\it #2}\endgraf
    \global\advance\conjno by 1}

\newcount\eqnno \eqnno=0
\def\eqn#1{\global\edef#1{(\the\secno.\the\eqnno)}#1\global
   \advance\eqnno by 1}
\newcount\secno \secno=0
\def\sec#1\par{\eqnno=1\global\advance\secno by 1\bigskip
    \noindent{\bf\the\secno. #1}\par\nobreak\noindent}

\def\definition{\noindent{\bf Definition.\enspace}}

\def\proof{\noindent{\bf Proof.\enspace}}

\def\qed{~\hbox{\quad\vbox{\hrule height.4pt
    \hbox{\vrule width.4pt height7pt \kern7pt \vrule width.4pt}
    \hrule height.4pt}}}

\parskip=11pt plus3pt minus3pt

\def\itc#1{{\it #1\/}}

\centerline{\bf What is a closed-form number?}
\bigskip
\centerline{Timothy Y. Chow}
\smallskip
\centerline{University of Michigan, Ann Arbor}

\sec Introduction

When I was a high-school student,
I liked giving \itc{exact} answers to numerical problems whenever possible.
If the answer to a problem were $2/7$ or $\pi\sqrt{5}$
or $\arctan 3$ or $e^{1/e}$, I would always leave it in that form
instead of giving a decimal approximation.

Certain problems frustrated me because there did not seem to be
any way to express their solutions exactly.
For example, consider the following problems.

\noindent
{\bf Question 1.} The equation
$$\eqalignno{x+e^x &= 0&\eqn\xex\cr}$$
has exactly one real root; call it~$R$.
Is there a closed-form expression for~$R$?

\noindent
{\bf Question 2.} The equation
$$\eqalignno{2x^5-10x+5 &= 0&\eqn\poly\cr}$$
has five distinct roots $r_1$, $r_2$, $r_3$, $r_4$, and~$r_5$.
Are there closed-form expressions for them?

\noindent
Questions like this seemed to have a negative answer,
but I continued hoping that the answer was yes,
and that I just did not know enough mathematics yet.

In college I learned about Galois theory,
and that the Galois group of equation~\poly\
is~$S_5$ [7,~Section~5.8].
So the $r_i$ are provably not expressible in terms of radicals.
But although this probably should have satisfied me, it did not.
Consider the equation
$$x^4-(6\sqrt3)x^3 + 8x^2 + (2\sqrt3)x - 1 = 0.$$
Its roots are $\tan(\pi/15)$, $\tan(4\pi/15)$, $\tan(7\pi/15)$,
and $\tan(13\pi/15)$.
These seemed to me to be perfectly good closed-form expressions.
Although in this particular case the roots could also be expressed
in terms of radicals, it seemed to me that there might exist
algebraic numbers that were not expressible using radicals
but that could still be expressed in closed form---say,
using trigonometric or exponential or logarithmic functions.
So as far as I was concerned, Galois theory was not the end of the story.

When students ask for a closed-form expression for $\int \exp(x^2)\,dx$,
we all know the standard answer:
the given function is not an elementary function.
Curiously, though, Question~1
(as well as Question~2,
if you accept my dissatisfaction with the Galois-theoretic answer)
does not seem to have a standard answer that ``everybody knows.''
At most we might mutter vaguely that equation~\xex\ is
a ``transcendental equation,'' but this is not very helpful.

This nonexistence of a standard answer to such a simple
and common question seems almost scandalous to me.
The main purpose of this paper is to eliminate this scandal by suggesting
a precise definition of a ``closed-form expression for a number.''
This will enable us to restate Questions 1 and~2 precisely,
and will let us see how they are related to existing work in
logic, computer algebra, and transcendental number theory.
My hope is that this definition of
a closed-form expression for a number will become standard,
and that many readers will be lured into working on
the many attractive open problems in this~area.

\sec From elementary functions to EL numbers

How can we make Questions 1 and~2 precise?
Our first inclination might be to turn to
the notion of an \itc{elementary function}.
Recall that a function is \itc{elementary} if it can be
constructed using only a finite combination of
constant functions, field operations,
and algebraic, exponential, and logarithmic functions.
This class of functions has been studied a great deal
in connection with the problem of symbolic integration
or ``integration in finite terms''~[4],
and it does a rather good job
of capturing ``high-school intuitions'' about
what a closed-form expression should look like.
For example, in Question~1 above, it turns out that $R = -W(1)$,
where $W$, the \itc{Lambert~$W$ function}~[6],
is the (multivalued) function defined by the equation
$$W(x)e^{W(x)} = x.$$
But since $W$ is not an elementary function~[5],
this is not an answer that would satisfy most high-school students.
Similarly, if we allow various special functions---e.g.,
elliptic, hypergeometric, or theta functions---then we can explicitly
express the $r_i$ in Question~2,
or indeed the roots of any polynomial equation,
in terms of the coefficients~[3, 10].
But this again feels unsatisfactory because these special functions are
not elementary.

The concept of an elementary function is certainly on the right track,
but observe that what we need for Questions 1 and~2
is a notion of a closed-form \itc{number}
rather than a closed-form \itc{function}.
The distinction is important;
we cannot, for example, simply define an ``elementary number''
to be any number obtainable by evaluating an elementary function
at a point, because all constant functions are elementary,
and this definition would make all numbers elementary.
Furthermore, even if a function (like~$W$) is not elementary,
it is conceivable that each particular value that it takes
($W(1), W(2), \ldots$) could have an elementary expression,
but with different-looking expressions at different points.
So we will not try to define closed-form numbers \itc{in terms of}
elementary functions, but will give an \itc{analogous} definition.

We mention one more technical point.
By convention, all algebraic functions are considered to be elementary,
but this is not suitable for our purposes.
Intuitively, ``closed-form'' implies ``explicit,''
and most algebraic functions have no simple explicit expression.
So the set of \itc{purely transcendental elementary functions}
is a better prototype for our purposes than the set of elementary functions.
(``Purely transcendental'' simply means that
the word ``algebraic'' is dropped from the definition.)

With all these considerations in mind,
we propose the following fundamental definition.

\definition
A subfield $F$ of~$\C$ is \itc{closed under $\exp$ and $\log$} if 
(1)~$\exp(x)\in F$ for all $x\in F$ and (2)~$\log(x)\in F$
for all nonzero $x\in F$,
where $\log$ is the branch of the natural logarithm function
such that $-\pi<{\rm Im}(\log x) \le \pi$ for all~$x$.
The field~$\E$ of \itc{EL numbers} is the intersection
of all subfields of~$\C$ that are closed under $\exp$ and $\log$.

Before discussing~$\E$, let us make some remarks about terminology.
It might seem more natural to call~$\E$ the field of \itc{elementary
numbers,} but unfortunately this term is already taken.
It seems to have been first used by Ritt~[18, p.~60].
By analogy with elementary functions,
Ritt thought of elementary numbers as the smallest
\itc{algebraically closed} subfield~$\L$ of~$\C$
that is closed under $\exp$ and $\log$.
It so happens that terminology has evolved since Ritt,
so that $\L$ is now known as the field of \itc{Liouvillian numbers,}
and ``elementary numbers'' are now numbers that can be specified
\itc{implicitly} as well as explicitly by exponential, logarithmic,
and algebraic operations~[16].
But either way, calling~$\E$ the field of elementary numbers
would conflict with existing usage.
The ``EL'' in the term ``EL number'' is intended to be an abbreviation for
``Exponential-Logarithmic'' as well as a diminutive of ``ELementary,''
reminding us that $\E$ is a subfield of the elementary numbers.

I should also remark that I am certainly not the first person ever
to have considered the field~$\E$,
but it has received surprisingly little attention in the literature
and nobody seems to have lobbied for it as a fundamental object of interest,
which in my opinion it is
(as illustrated by my temerity in using ``blackboard bold'' for~it).

Let us do a few warmup exercises to familiarize ourselves with~$\E$.
We can construct $\E$ as follows.  Set $\E_0=\{0\}$, and for each $n>0$
let $\E_n$ be the set of all complex numbers obtained
either by applying a field operation to any pair of
(not necessarily distinct) elements of~$\E_{n-1}$
or by applying $\exp$ or $\log$ to any element of~$\E_{n-1}$.
(Of course, division by zero and taking the logarithm of zero
are forbidden.)
Then it is clear that $\E$ is the union of all the~$\E_n$.
This shows in particular that $\E$ is countable,
and that every element of~$\E$ admits an explicit
finite expression in terms of rational numbers, field operations,
$\exp$, and~$\log$.

Most familiar constants lie in~$\E$, e.g.,
$$\eqalign{e &= \exp(\exp(0))\cr
   i &= \exp\biggl( {\log(-1) \over 2} \biggr)\cr
 \pi &= -i\log(-1)\cr}$$
Since $2\pi i \in \E$, we actually have access to all
branches of the logarithm and not just the principal one,
so all~$n$ of the $n$th roots of any $x\in\E$ are also in~$\E$.
It follows that all the roots of any polynomial
equation with rational coefficients that is solvable in radicals
lie in~$\E$.
Finally, formulas like
$$\eqalignno{ x^{2/3} &= \exp\biggl({2 \log x \over 3}\biggr) \cr
  \sin x &= {\exp(ix) - \exp(-ix) \over 2i} \cr
  \tanh x &= {\exp(x) - \exp(-x) \over \exp(x) + \exp(-x)}\cr
  \arccos x &= -i \log\biggl(x +
          \exp\biggl({\log(x^2-1)\over 2}\biggr)\biggr) \cr}$$
show that any expression involving ``high-school'' functions
and elements of~$\E$ is also in~$\E$.

We hope that this brief discussion has persuaded the reader
that $\E$ is the ``right'' precise definition of
``the set of all complex numbers that can be written in closed form.''
Accepting this, we can reformulate Questions 1 and~2 as follows.

\conjecture\conjone
{The real root~$R$ of $x+e^x = 0$ is not in~$\E$.}

\conjecture\conjtwo
{The roots $r_1$, $r_2$, $r_3$, $r_4$, and $r_5$
of $2x^5-10x+5 = 0$ are not in~$\E$.}

As far as I know, \conjone\ and \conjtwo\
are---perhaps surprisingly---still open.
Thus we are still frustrated,
but at least our frustration has been raised to a higher plane.
The next section of this paper is devoted to partial results.

\sec Schanuel's conjecture

\conjone\ is essentially due to Ritt,
except that he asked the question with $\L$ instead of~$\E$,
since he was motivated by different considerations from ours.
The best partial result I am aware of is due to Ferng-Ching Lin~[12].
To state Lin's theorem, we must first recall \itc{Schanuel's conjecture}.

\noindent
{\bf Schanuel's Conjecture}.
{\it
If $\alpha_1, \alpha_2, \ldots, \alpha_n$ are complex numbers
linearly independent over~$\Q$, then the transcendence
degree of the field
$\Q\bigl(\alpha_1, e^{\alpha_1}, \alpha_2, e^{\alpha_2}, \ldots,
 \alpha_n, e^{\alpha_n}\bigr)$
over~$\Q$ is at least~$n$.}

To orient the reader who has not seen Schanuel's conjecture before,
we mention that it implies many famous theorems and
conjectures about transcendental numbers.
For example, Schanuel's conjecture implies
the Lindemann-Weierstrass theorem [2, Theorem~1.4] that if
$\alpha_1, \alpha_2, \ldots, \alpha_n$ are algebraic numbers
that are linearly independent over~$\Q$,
then $e^{\alpha_1}, e^{\alpha_2}, \ldots, e^{\alpha_n}$
are algebraically independent over~$\Q$.
Schanuel's conjecture also implies the Gelfond-Schneider theorem
that if $\alpha_1$ and~$\alpha_2$ are algebraic numbers
for which there exist $\Q$-linearly independent
numbers $\beta_1$ and~$\beta_2$
such that $\alpha_1 = e^{\beta_1}$ and $\alpha_2 = e^{\beta_2}$,
then $\beta_1$ and~$\beta_2$ are linearly independent over
the algebraic numbers.
Baker's generalization~[2, Theorem~2.1] of Gelfond-Schneider to
an arbitrarily large finite number of~$\alpha_i$
also follows from Schanuel's conjecture.
It is an easy exercise (using $e^{\pi i} = -1$)
to show that Schanuel's conjecture implies that $e$ and~$\pi$
are algebraically independent, which is currently not known.
(It is not even known that $e+\pi$ is transcendental.)
A proof of Schanuel's conjecture would be big news,
although at present it seems to be out of reach.

Let $\overline \Q$ denote the algebraic closure of~$\Q$.
Then Lin's result is the following.

\theorem\linthm
{If Schanuel's conjecture is true and $f(x,y) \in \overline\Q[x,y]$
is an irreducible polynomial involving both $x$ and~$y$ and
$f\bigl(\alpha, \exp(\alpha)\bigr) = 0$ for some nonzero $\alpha\in\C$,
then $\alpha\notin\L$.}

By taking $f(x)=x+y$ and noting that $\E\subseteq\L$,
we see at once that Schanuel's conjecture implies \conjone.

\conjtwo\ seems to be new.
The literature does contain some negative results about
the insolubility of \itc{general} polynomial equations
in terms of the exponential and logarithmic functions,
although it is difficult to find a satisfactory reference
(the best I know are [9, paragraph~513] and [1, p.~114]).
The inexpressibility of an algebraic function
in terms of $\exp$ and $\log$ does not, however,
imply that \itc{particular values} of an algebraic function
cannot be expressed in terms of $\exp$ and $\log$,
just as some quintic equations with rational coefficients
are solvable in radicals even though the general quintic is not.
Since \conjtwo\ has not been seriously attacked before,
it may turn out to be quite accessible.

The remainder of this section is devoted to proving the following result.

\theorem\mainthm
{Schanuel's conjecture implies \conjone\ and \conjtwo.}

As we just remarked, Lin has already shown that Schanuel's conjecture
implies \conjone, but we shall exploit the fact that \conjone\ is
weaker than the conclusion of \linthm\ to give a shorter proof.
The proof of \mainthm\ below is joint work with Daniel Richardson.
Inspection of the proof shows that it is readily generalized to
other transcendental equations that are similar to~\xex\
and to any algebraic number whose minimum polynomial has
a Galois group that is not solvable.
In particular, if Schanuel's conjecture is true,
then our notion of a ``closed-form algebraic number''
coincides with the usual one, i.e., solvability in radicals.

Although \conjone\ and \conjtwo\ involve
quite different kinds of equations,
it turns out that there is a single concept
(that of a \itc{reduced tower,} defined below) that is the key to both.
The reader who is not interested in the details of the proof
may skip directly to the next section now without loss of continuity.

We need some preliminaries.
If $A = (\alpha_1, \alpha_2, \ldots, \alpha_n)$
is a finite sequence of complex numbers
then for brevity we write $A_i$ for the field
$\Q(\alpha_1, e^{\alpha_1}, \alpha_2, e^{\alpha_2},
\ldots, \alpha_i, e^{\alpha_i})$.
In particular, $A_0 = \Q$.

\definition
A \itc{tower} is a finite sequence
$A = (\alpha_1, \alpha_2, \ldots, \alpha_n)$
of nonzero complex numbers such that for all~$i \in \{1, 2, \ldots, n\}$,
there exists some integer $m_i>0$
such that $\alpha_i^{m_i} \in A_{i-1}$ or
$e^{\alpha_i m_i} \in A_{i-1}$ (or both).
A tower is \itc{reduced} if the set $\{\alpha_i\}$ is linearly
independent over~$\Q$.
If $\beta \in \C$, then a \itc{tower for~$\beta$}
is a tower $A = (\alpha_1, \alpha_2, \ldots, \alpha_n)$
such that $\beta\in A_n$.

For any $\gamma\in\E$, there exists a tower for~$\gamma$.
This is best explained by example.
Suppose $\gamma=4 + \log\bigl(1 + e^{(\log 2)/3}\bigr)$.
Then we may take
$$A = (\alpha_1, \alpha_2, \alpha_3)
    = \bigl(\log 2, (\log 2)/3, \log(1 + e^{(\log 2)/3})\bigr).$$
We can then take $m_i = 1$ for all~$i$, because
$e^{\alpha_1} = 2 \in A_0$, $\alpha_2 \in A_1$, and $e^{\alpha_3}\in A_2$.
In general, we build up the expression for~$\gamma$ step by step,
and if at step~$i$ we need to take the exponential of some
number~$\beta \in A_{i-1}$ we simply set $\alpha_i=\beta$;
if we need to take the logarithm of some $\beta\in A_{i-1}$
then we set $\alpha_i = \log \beta$.
With this construction, we never need to take $m_i > 1$,
but the tower we obtain may not be reduced
(as is the case in this example: $\alpha_1 - 3\alpha_2 = 0$).
In order to be able to use Schanuel's conjecture,
however, we need reduced towers,
so our first goal is to show how to reduce a given tower.

\noindent
{\bf Division Lemma.}
{\it
Suppose $A = (\alpha_1, \alpha_2, \ldots, \alpha_n)$
is a tower and $q_1, q_2, \ldots, q_n$
are nonzero integers.  Then the sequence
$B = (\beta_1, \beta_2, \ldots, \beta_n)$ defined by
$\beta_i = \alpha_i/q_i$ is also a tower,
and $A_i \subseteq B_i$ for all~$i$.}

\proof
Given any $i \in \{1, 2, \ldots, n\}$,
note that every $\gamma\in A_i$ is a rational function
(with rational coefficients)
of the numbers $\alpha_1, e^{\alpha_1}, \ldots, \alpha_i, e^{\alpha_i}$.
Now
$$\alpha_j = (\alpha_j/q_j)q_j = \beta_j q_j \qquad {\rm and} \qquad
  e^{\alpha_j} = e^{(\alpha_j/q_j)q_j} =
   (e^{\beta_j})^{q_j} \qquad \hbox{for all~$j$,}$$
so $\gamma$ is also a rational function with rational coefficients
of the numbers $\beta_1, e^{\beta_1}, \ldots, \beta_i, e^{\beta_i}$,
and hence $\gamma\in B_i$.  So $A_i \subseteq B_i$ for all~$i$.

Given any $i \in \{1, 2, \ldots, n\}$,
there is some integer $m_i > 0$ such that $\alpha_i^{m_i} \in A_{i-1}$
or $e^{\alpha_i m_i} \in A_{i-1}$.
Consider first the case in which $\alpha_i^{m_i} \in A_{i-1}$.  Then
$$\beta_i^{m_i} = \biggl({\alpha_i \over q_i}\biggr)^{m_i} \in A_{i-1}
   \subseteq B_{i-1}.$$
If on the other hand $e^{\alpha_i m_i} \in A_{i-1}$, then
$$e^{\beta_i (q_i m_i)} = e^{\alpha_i m_i} \in A_{i-1} \subseteq B_{i-1}.$$
Hence there is a positive integer $m'_i$ (for example, $m'_i = q_i m_i$)
such that $\beta_i^{m'_i} \in B_{i-1}$ or $e^{\beta_i m'_i} \in B_{i-1}$.
So $B$ is a tower.\qed

\vbox{
\noindent
{\bf Reduction Lemma.}
{\it For any $\gamma\in\E$, there exists a reduced tower for~$\gamma$.}

\proof
If $\gamma\in\Q$ then we may take $A$ to be the empty sequence.
Otherwise, suppose that every tower for~$\gamma$ is not reduced;
we shall derive a contradiction.
Choose such an~$A$ with $n$ minimal; since $\gamma\notin\Q$, $n\ge1$.
Let $i$ be the smallest integer such that
$\{\alpha_1, \alpha_2, \ldots, \alpha_i\}$ is linearly dependent.
Then
$$\eqalignno{\alpha_i &= \sum_{j=1}^{i-1}
   {p_j \alpha_j \over q_j}&\eqn\reduct\cr}$$
for some integers $p_1, q_1, p_2, q_2, \ldots, p_{i-1}, q_{i-1}$.
We claim that the sequence
$$A' = \biggl( {\alpha_1 \over q_1}, {\alpha_2 \over q_2}, \ldots,
  {\alpha_{i-1} \over q_{i-1}}, \alpha_{i+1}, \alpha_{i+2},
     \ldots, \alpha_n\biggr)$$
is a tower for~$\gamma$.  Since $A'$ is shorter than~$A$,
this contradicts the minimality of~$n$ and proves the theorem.
}

To prove the claim, note first that by the division lemma, the sequence
$$\biggl( {\alpha_1 \over q_1}, {\alpha_2 \over q_2}, \ldots,
  {\alpha_{i-1} \over q_{i-1}}\biggr)$$
is a tower.  Next, note that equation~\reduct\ implies that
$\alpha_i \in A'_{i-1}$ and also, by exponentiating, that
$e^{\alpha_i}$ is a polynomial (in fact a monomial) in the numbers
$e^{\alpha_1/q_1}, \ldots, e^{\alpha_{i-1}/q_{i-1}}$, so that
$e^{\alpha_i} \in A'_{i-1}$.
By the division lemma, $A_{i-1} \subseteq A'_{i-1}$, so
$$A'_{i-1} \supseteq A_{i-1}(\alpha_i, e^{\alpha_i}) = A_i.$$
This ensures that the tower condition for~$A'$ is satisfied
at the boundary between $\alpha_{i-1}/q_{i-1}$ and~$\alpha_{i+1}$,
and also that $A'_{n-1} \supseteq A_n \ni \gamma$,
proving the claim.\qed

\noindent
{\bf Proof of \mainthm.}
We first make a general remark.
If $B = (\beta_1, \beta_2, \ldots, \beta_n)$ is a reduced tower,
then Schanuel's conjecture implies that for all~$i$,
\itc{exactly one} of $\beta_i$ and~$e^{\beta_i}$
is algebraic over~$B_{i-1}$.
For by the definition of a tower, \itc{at least one}
of the two is algebraic over~$B_{i-1}$;
this implies that the transcendence degree of~$B_i$ over~$\Q$
is at most~$i$ for all~$i$.
Then because $B$ is reduced,
Schanuel's conjecture may be applied,
telling us that $\alpha_i$ and $e^{\alpha_i}$ cannot \itc{both}
be algebraic over~$B_{i-1}$,
and that the transcendence degree of~$B_i$ over~$\Q$
is exactly~$i$.

Now assume Schanuel's conjecture.
We first prove \conjone.
Assume $R\in\E$; we derive a contradiction.
By the reduction lemma, there is a reduced tower
$A=(\alpha_1, \alpha_2, \ldots, \alpha_n)$ for~$R$.
Since $e$ is transcendental, $R\notin\Q$
(for $R = p/q$ implies $e^p = (-p/q)^q$), so $n\ge1$.
By truncating the tower if necessary,
we may assume that $R\notin A_i$ if $i<n$.

Let $A' = (\alpha_1, \alpha_2, \ldots, \alpha_n, R)$.
Then $R\in A'_n$, 
and the relation $R + e^R = 0$
shows that $e^R \in A'_n$ as well.
By our ``general remark,'' $A'$ cannot be reduced.
But $A$ is reduced, so
$$R = \sum_{i=1}^n {p_i \alpha_i \over q_i}$$
for some integers $p_1, q_1, p_2, q_2, \ldots, p_n, q_n$.
Moreover, $p_n\ne0$ because $R \notin A_i$ for $i<n$.
The relation $R + e^R = 0$ becomes
$$\eqalignno{ \sum_{i=1}^n {p_i \alpha_i \over q_i}
  + \prod_{i=1}^n \bigl(e^{\alpha_i/q_i}\bigr)^{p_i} &= 0.&\eqn\gpluseg\cr}$$
Let $A'' = (\alpha_1/q_1, \alpha_2/q_2, \ldots, \alpha_n/q_n)$.
By the division lemma, $A''$ is a tower,
and since $A$ is reduced, $A''$ is reduced.
But since $p_n\ne 0$, equation~\gpluseg\
shows that if $\alpha_n/q_n$ is algebraic over~$A''_{n-1}$
then so is $e^{\alpha_n/q_n}$, and vice versa.
By our ``general remark,'' $A''$ cannot be reduced,
and this gives our desired contradiction.

Now for \conjtwo.
We shall assume that the reader is familiar with
the rudiments of Galois theory.
Assume that $r_1 \in \E$; we derive a contradiction.
By the reduction lemma there is a reduced tower
$A=(\alpha_1, \alpha_2, \ldots, \alpha_n)$ for~$r_1$
(of course, unrelated to the ``$A$'' in the first part of this proof).
For all~$i$, let
$$\beta_i =
\cases{\alpha_i,&if $\alpha_i$ is transcendental over $A_{i-1}$;\cr
    e^{\alpha_i},&if $e^{\alpha_i}$ is transcendental over $A_{i-1}$.\cr}$$
Then the $\beta_i$ are algebraically independent and form a transcendence
basis for $A_n$ over~$\Q$.  Let $F=\Q(\beta_1, \beta_2, \ldots, \beta_n)$.
Clearly $A_n$ is an extension by radicals of~$F$.
Let $L$ be the Galois closure of~$A_n$ over~$F$;
then ${\rm Gal}(L/F)$ is solvable.
If $F'$ is the splitting field over~$\Q$ of the polynomial~\poly,
then ${\rm Gal}(F'\!/\Q) = S_5$
and $F'\cap F = \Q$ because $F'\!/\Q$ is an algebraic extension
while $F/\Q$ is a purely transcendental extension.
Therefore the compositum $FF'$ is Galois over~$F$ with
Galois group~$S_5$ [11, Chapter~8, Theorem~1.12].
But $FF' \subseteq L$, so $S_5$ must be a homomorphic image of
${\rm Gal}(L/F)$.  This is our desired contradiction,
because every homomorphic image of a solvable group is solvable.\qed

\sec Related work and open problems

Much of the work that has been done on fields like $\E$, $\L$,
or the field of elementary numbers has been motivated by
problems in logic and computer algebra.
A typical problem is this:
given a complicated expression for a number in~$\E$,
how can you tell if it exactly equals zero?
Clearly this is an important problem for designers
of symbolic computation software.
It is harder than it might seem at first glance,
and is still not fully solved, although
Richardson~[17] has explicitly described a procedure that 
takes a given elementary number and, if the procedure terminates,
correctly says whether or not the number equals zero.
He has also proved that if Schanuel's conjecture is true,
then the procedure does in fact always terminate.
This more or less solves the zero-recognition problem for elementary numbers
(and \itc{a~fortiori} for $\E$ and~$\L$) in practice.

The zero-recognition problem is closely related to
a famous long-standing question of Tarski.
Tarski proved that the first-order theory of the real numbers is decidable,
which implies in particular that there is an algorithm
for determining whether or not
any given finite system of polynomial equations and inequalities
has a solution in the reals~[8, p.~340].
The proof proceeds by \itc{quantifier elimination,}
which we can think of roughly as follows:
the statement that ``there exists a solution''
involves existential quantifiers,
and quantifier elimination is a procedure for
transforming such statements into ones that are quantifier-free.
These are then easy to check because
all that is involved is a zero-recognition problem for integers.
After proving his theorem,
Tarski asked if it could be extended to
the first-order theory of the real numbers with exponentiation.

This problem is very hard, because it turns out that
quantifier elimination is not possible in this theory.
Moreover, checking quantifier-free statements involves
the zero-recognition problem for expressions with exponentials,
which is difficult.
Great progress has recently been made, however.
Macintyre~[13] showed that if Schanuel's conjecture is true,
then there is a decision procedure for the quantifier-free statements.
Then Wilkie~[19] proved in~1991 that
the first-order theory of reals with exponentiation is model complete
(which roughly means that quantifiers can ``almost'' be entirely eliminated).
Building on this work, Macintyre and Wilkie~[14]
showed that if Schanuel's conjecture is true, then
the first-order theory of the real numbers with exponentiation is decidable.  
In particular, from these methods one can extract a zero-recognition
procedure for elementary numbers (again, contingent on Schanuel).
See~[15] for a splendid account of these and related results.

Zero-recognition in~$\E$ should be easier than
zero-recognition in the elementary numbers.
Can one recognize zero in~$\E$ without assuming Schanuel,
or at least by assuming something weaker?
The ideas of Macintyre~[13] are a good starting-point here.

Another interesting open problem,
posed by Thomas Colthurst (in a {\tt sci.math} article
posted on June 21, 1993),
is to produce an explicit example of a number that is not in~$\E$.
Since $\E$ is countable, Cantor's diagonal argument
gives us an algorithm for producing the decimal expansion
of a non-EL number, but this is not very satisfying.
Colthurst suggests that one might expect an expression of the form
$$F = \sum_{m=1}^\infty f(m)$$
to work, where $f(m)$ is a nonnegative function of~$m$
that approaches zero rapidly.
For example, the sets $\E_n$ (defined above in section~2) are finite,
so there exists $\epsilon_n > 0$ such that any two distinct
numbers in~$\E_n$ differ by at least~$\epsilon_n$.
If one could find $f$ with the property that, for all~$n$,
$$\sum_{m=1}^n f(m) \in \E_n \qquad{\rm and}\qquad
  \sum_{m=n+1}^\infty f(m) < \epsilon_n,$$
then $F$ could not be in~$\E_n$ for any~$n$.
This is probably too na\"\i ve, but it seems that something
along these lines should be feasible.  Can $f$ be chosen
to be elementary?

On the grounds that many high-school students are unfamiliar with
complex numbers, one can ask for a ``real analogue'' of~$\E$.
What is the right definition?  Such a real analogue would lack
many of the nice properties of~$\E$ (e.g., recall that if
an irreducible cubic with rational coefficients
has three distinct real roots,
then they cannot be expressed using radicals alone
if complex numbers are forbidden), but it might still be interesting.

Finally, we mention that Richardson (personal communication)
has shown that if Schanuel's conjecture is false,
then there is a counterexample involving only elementary numbers.
Can this be strengthened to show that any counterexample
must lie in~$\E$?

We hope the reader will be tempted to attack
these relatively untouched questions.

\sec Acknowledgments

The notion of an EL number occurred to me a long time ago and I have
benefited from discussions with numerous people over the years.
The ideas of Daniel Richardson and Thomas Colthurst
have been particularly helpful and have had a profound
influence on this paper.
I would like to thank 
Alexander Barvinok, Jerry Shurman, Arpad Toth, and
Robert Corless for helpful comments on the mathematics
related to references [1, 5, 6, 10].
Thanks also to
Juergen Weiss, Alexander Pruss, and David Feldman
for brief email messages that they may have forgotten by now but
which helped point me in the right direction when I was getting started.
Finally, this work was supported in part by
a National Science Foundation Postdoctoral Fellowship.

\sec References

\item{1.}
V. B. Alekseev, \itc{Abel's Theorem in Problems and Solutions,}
Izdat.~``Nauka,'' 1976 (Russian).

\item{2.}
A. Baker, \itc{Transcendental Number Theory,}
Cambridge Mathematical Library,
Camb.\ Univ.\ Press, 1990.

\item{3.}
G. Belardinelli, Fonctions hyperg\'eom\'etriques de plusieurs
variables et r\'esolution analytique des \'equations alg\'ebriques
g\'en\'erales, \itc{M\'emorial des Sci.\ Math.}\ {\bf 145} (1960).

\item{4.}
M. Bronstein, \itc{Symbolic Integration I: Transcendental Functions,}
Algorithms and Computation in Mathematics, Volume~1, Springer-Verlag, 1997.

{
\catcode`\_=12
\font\smtt=cmtt9

\item{5.}
R. M. Corless, Is elementary? Math 498/990 notes, Nov.~23, 1995,
URL\hfill\break
{\smtt http://www.apmaths.uwo.ca/\~{ }rmc/AM563/NOTES/Nov_23_95/Nov_23_95.html
}

}

\item{6.}
R. M. Corless, G. H. Gonnet, D. E. G. Hare, D. J. Jeffrey, and
D. E. Knuth, On the Lambert $W$ function, \itc{Adv.\ Comput.\ Math.}\
{\bf 5} (1996), 329--359.


\item{7.}
I. N. Herstein, \itc{Topics in Algebra,} 2nd ed., Wiley, 1975.

\item{8.}
N. Jacobson, \itc{Basic Algebra I,} 2nd ed., W. H. Freeman, 1985.

\item{9.}
C. Jordan, \itc{Trait\'e des Substitutions et des \'Equations
Alg\'ebriques,} Gauthier-Villars, 1870.

\item{10.}
R. Bruce King, \itc{Beyond the Quartic Equation,} Birkh\"auser Boston, 1996.

\item{11.}
S. Lang, \itc{Algebra,} 2nd ed., Addison-Wesley, 1984.

\item{12.}
F.-C. Lin, Schanuel's conjecture implies Ritt's conjecture,
\itc{Chinese J. Math.}\ {\bf 11} (1983), 41--50.

\item{13.}
A. Macintyre, Schanuel's conjecture and free exponential rings,
\itc{Ann.\ Pure Appl.\ Logic} {\bf 51} (1991), 241--246.

\item{14.}
A. Macintyre and A. J. Wilkie,
On the decidability of the real exponential field,
in \itc{Kreiseliana: About and Around Georg Kreisel,}
ed.~P.~Odifreddi, A. K. Peters, 1996.

\item{15.}
D. Marker, Model theory and exponentiation,
\itc{Notices Amer.\ Math.\ Soc.}\ {\bf 43} (1996), 753--759.

\item{16.}
D. Richardson, The elementary constant problem, in \itc{Proceedings
of the International Symposium on Symbolic and Algebraic Computation,
Berkeley, July 27--29, 1992,} ed.~P. S. Wang, ACM Press, 1992.

\item{17.}
D. Richardson, How to recognize zero, \itc{J. Symb.\ Comp.}\
{\bf 24} (1997), 627--645.

\item{18.}
J. Ritt, \itc{Integration in Finite Terms: Liouville's Theory of Elementary
Models,} Colum\-bia Univ.\ Press, 1948.

\item{19.}
A. J. Wilkie, Model completeness results for expansions of the ordered field
of real numbers by restricted Pfaffian functions and the exponential
function, \itc{J.~Amer.\ Math.\ Soc.}\ {\bf 9} (1996), 1051--1094.
\bye